\documentclass[a4paper,11pt]{article}

\usepackage{amssymb,amsmath,amsthm}
\usepackage{setspace}
\usepackage{graphicx} 
\usepackage{enumerate}
\usepackage{chngcntr}
\usepackage{mathtools}
\counterwithin{equation}{section}
\usepackage{marginnote}
\usepackage{pgf,tikz}
\usepackage{mathrsfs}
\usetikzlibrary{arrows}
\usetikzlibrary[patterns]
\tikzstyle{dashed}=[dash pattern=on 3pt off 3pt]

\usepackage{tikz-cd}
\usetikzlibrary{%
  matrix,%
  calc,%
  arrows%
}

\usepackage{verbatim}

\usepackage[backend=bibtex,maxbibnames=99]{biblatex}

\usepackage{filecontents}
\usepackage{sectsty}
\usepackage{lipsum}

\theoremstyle{definition}

\newtheorem{thm}{Theorem}[section]
\newtheorem{prop}[thm]{Proposition}
\newtheorem{coro}[thm]{Corollary}
\newtheorem{lem}[thm]{Lemma}
\newtheorem{dfn}[thm]{Definition}

\newtheorem{con}[thm]{Condition}
\newtheorem{ex}[thm]{Example}

\newtheorem{rmk}[thm]{Remark}

\allsectionsfont{\centering}

\providecommand{\lan}{\mathcal{L}}

\providecommand{\M}{\mathcal{M}}
\providecommand{\cl}{\mathcal{C}}
\providecommand{\rat}{\mathbb{Q}}

\newcommand{\rel}{\mathbf{r}}

\def\Ind#1#2{#1\setbox0=\hbox{$#1x$}\kern\wd0\hbox to 0pt{\hss$#1\mid$\hss}
\lower.9\ht0\hbox to 0pt{\hss$#1\smile$\hss}\kern\wd0}
\def\Notind#1#2{#1\setbox0=\hbox{$#1x$}\kern\wd0\hbox to 0pt{\mathchardef
\nn="3236\hss$#1\nn$\kern1.4\wd0\hss}\hbox to 0pt{\hss$#1\mid$\hss}\lower.9\ht0
\hbox to 0pt{\hss$#1\smile$\hss}\kern\wd0}
\def\ind{\mathop{\mathpalette\Ind{}}}

\begin{document}

\title{Automorphism groups of homogeneous structures with stationary weak independence relations\footnote{This is an ongoing work. Further work that is being done is discussed in Section 4.}}

\author{Yibei Li}

\newcommand{\Addresses}{{
  \bigskip
  \footnotesize

  Yibei Li, \textsc{Department of Mathematics, Imperial College London,
    Huxley Building, 180 Queens Gate, London SW7 2AZ}\par\nopagebreak
  \textit{E-mail address}: \texttt{yibei.li16@imperial.ac.uk}

}}

\date{}

\maketitle

\section*{Abstract}

We generalise the result of \cite{tentziegler2013isometry} to homogeneous structures that have a stationary independence relation without the symmetry axiom. We apply our result to prove simplicity of the automorphism group of some asymmetric examples due to Cherlin \cite{cherlin1998classification}.

\section{Introduction}

Given a relational language $\lan$, a countable $\lan$-structure $\M$ is \emph{homogeneous} if every partial isomorphism between finite substructures of $\M$ extends to an automorphism of $\M$. Fra{\"\i}ss{\'e}'s Theorem \cite{fraisse1953theorem} provides one way of constructing homogeneous structures by establishing a one-to-one correspondence between such structures and amalgamation classes. We call the homogeneous structure the \emph{Fra{\"\i}ss{\'e} limit} of the corresponding amalgamation class.


Given finite $\lan$-structures $A,B,C$ where $B \subseteq A,C$, the \emph{free amalgam} of $A,C$ over $B$ is the $\lan$-structure $D$ on the disjoint union of $A,C$ over $B$ and for each relation $R \in \lan$, $R^D=R^A \cup R^C$. An amalgamation class $\cl$ is \emph{free} if it is closed under taking free amalgams. A homogeneous structure is \emph{free} if it is the Fra{\"\i}ss{\'e} limit of a free amalgamation class. In \cite{macphersontent2011simplicity}, Macpherson and Tent proved the following theorem about free homogeneous structures using ideas and results from model theory and topological groups:

\begin{thm}\label{tztheorem}(\cite{macphersontent2011simplicity})
Let $\M$ be a countable free homogeneous relational structure. Suppose $Aut(\M) \neq Sym(\M)$ and $Aut(\M)$ is transitive on $\M$. Then $Aut(\M)$ is simple.
\end{thm}

This is then generalised by Tent and Ziegler \cite{tentziegler2013isometry} to a weaker notion than free homogeneous structures, namely a homogeneous structure with a \emph{stationary independence relation} (see Definition \ref{swir}). Tent and Ziegler then applied their method to the Urysohn space \cite{tentziegler2013isometry} as well as the bounded Urysohn space \cite{tentziegler2013bounded}, which are not free, but have local stationary independence relations.

In the appendix of \cite{cherlin1998classification}, Cherlin proposed a generalisation of free amalgamation class, called \emph{semi-free amalgamation class} (see Definition \ref{semifree}). He provided some examples for languages consisting either symmetric or asymmetric relations. It is then natural to ask whether Tent and Ziegler's method could be applied to semi-free homogeneous structures. The author of this paper showed in \cite{li2018simplicity} that their method can be applied to all the symmetric examples in the appendix of \cite{cherlin1998classification} as well as some general cases and proved the simplicity of the automorphism groups of the Fra{\"\i}ss{\'e} limit of those semi-free amalgamation class. 
The same statement was proved in \cite{evans2019simplicity} for more general structures using the method from \cite{tentziegler2013bounded}. However, we cannot apply the method to the asymmetric examples as they do not satisfy the symmetry axiom. The main result of this paper, stated below, is motivated by trying to generalise Tent and Ziegler's method to the asymmetric examples of Cherlin. 

We will generalise Tent and Ziegler's method to a stationary independence relation without the symmetry axiom in Section 2, which we define to be \emph{a stationary weak independence relation}. The main theorem we prove is the following.

\begin{thm}\label{main}
Suppose $\M$ is a countable structure with a stationary weak independence relation and $g \in Aut (\M)$ is such that $g$ and $g^{-1}$ move both almost $R$-maximally and almost $L$-maximally. Then any element of $G$ is a product of eight conjugates of $g$.
\end{thm}

We will then apply the method to some of the asymmetric examples in the appendix of \cite{cherlin1998classification} in Section 3 and show that their automorphism groups are simple. 

We now fix some notation for the paper. We first fix a first-order relational language $\lan$, which is specified by a set of relation symbols $\{R_i : i \in I \}$ and each $R_i$ has arity $r_i \in \mathbb{N}$. Then an $\lan$-structure is a set $A$ together with a subset $R^{A}_i \subseteq A^{r_i}$ for each $i \in I$ representing the structure on $A$. Let $\M$ be an $\lan$-structure and $A,B$ be finite substructures of $\M$, we use the notation $AB$ to denote the substructure of $\M$ on the underlying set $A \cup B$. We also simplify the notation $\{a\}B$ to $aB$. Let $G=Aut(\M)$ and denote the pointwise stabiliser of $B$ by $G_{(B)}$. For $g, h\in G$, let $g^h$ denote $hgh^{-1}$ and $[g,h]=g^{-1}h^{-1}gh$ be the commutator of $g,h$. For a homogeneous structure $\M$ and a finite subset $B\subseteq \M$, the model theoretic notion of an \emph{$n$-type over $B$} corresponds to a $G_{(B)}$-orbit of an $n$-tuple. For $a \in \M^n$, the \emph{type of $a$ over $B$}, denoted by $tp(a/B)$, is the type over $B$ whose corresponding $G_{(B)}$-orbit contains $a$. So, we may use $tp(a/B)$ to denote its corresponding $G_{(B)}$-orbit. Note that $a,a'$ have the same type over $B$ if they lie in the same $G_{(B)}$-orbit, i.e. there exists an automorphism of $\M$ that takes $a$ to $a'$ and fixes $B$ pointwise. We say $a$ \emph{realises} some type $p$ over $B$ if it lies in the corresponding $G_{(B)}$-orbit. We say a type is \emph{algebraic} if its set of realisations is finite, and \emph{non-algebraic} otherwise. In all of our examples, $tp(a/B)$ is algebraic if and only if $a \in B$.

\subsection*{Acknowledgement}
I would like to thank my supervisor, Prof. David Evans, for his supervision and guidance, which greatly assisted the research.

\section{Stationary Weak Independence Relation}

In this section, we follow a similar approach as in \cite{tentziegler2013isometry} to show that if $\M$ is a countable structure with \emph{a stationary weak independence relation}, defined as the following, then any element of $Aut(\M)$ is a product of conjugates of certain special automorphisms of $\M$.

\begin{dfn}\label{swir}
Let $\M$ be a homogeneous structure and suppose $A\ind _B C$ is a ternary relation between finite substructures $A,B,C$ of $\M$. We say that $\ind$ is a \emph{stationary weak independence relation} if the following axioms are statisfied:
\begin{enumerate}[(i)]
\item Invariance: for any $g \in Aut(\M)$, if $A \ind _B C$, then $ gA \ind _{gB} gC$

\item Monotonicity: $A \ind _B CD \Rightarrow A \ind _B C$, $A \ind _{BC} D$

\hspace{2.4cm} $AD \ind _B C \Rightarrow A \ind _B C$, $D \ind _{AB} D$

\item Transitivity: $A \ind _B C$, $A \ind _{BC} D \Rightarrow A \ind _B D $

\hspace{2.1cm} $A \ind _B C$, $D \ind _{AB} C \Rightarrow D \ind _B C $

\item Existence: If $p$ is an $n$-type over $B$ and $C$ is a finite set, then $p$ has a realisation $a$ such that $a \ind _B C$.

\hspace{1.8cm} If $p$ is an $n$-type over $B$ and $C$ is a finite set, then $p$ has a realisation $a$ such that $C \ind _B a$.

\item Stationarity: If $a$ and $a'$ are $n$-tuples that have the same type over $B$ and $a \ind _B C$, $a' \ind _B C$, then $a$ and $a'$ have the same type over $BC$.

\hspace{2.2cm} If $a$ and $a'$ are $n$-tuples that have the same type over $B$ and $C \ind _B a$, $C \ind _B a'$, then $a$ and $a'$ have the same type over $BC$.

\end{enumerate}

If in addition, $\M$ satisfies symmetry, i.e. $A\ind _B C  \Rightarrow C \ind _B A$, then we say $\ind$ is a stationary independence relation.

\end{dfn}

In \cite{kaplansimon2017automorphism}, Kaplan and Simon also generalised the notion of stationary independence relation, which they called a canonical independence relation. They studied the automorphism groups of structures with such an independence relation.

\begin{ex}\label{dloex}
Let $\M$ be the dense linear order $(\rat, \leq)$. For finite substructures $A,B,C \subseteq \rat$ such that $A \cap C \subseteq B$, we can define $A \ind _B C$ if for any $a\in A \setminus B, c \in C \setminus B$ such that $a \leq c$, then there exists $b\in B$ such that $a\leq b \leq c$. Then it can be checked that it is a stationary weak independence relation on $\rat$. 

For transitivity, let $B \subseteq A \subseteq \M, B\subseteq C \subseteq D \subseteq \M$. Suppose $A\ind _B C$ and $A\ind _C D$. We want to show $A \ind _B D$. For $a\in A \setminus C, d \in D \setminus C$ such that $a \leq d$, then there exists $c\in C$ such that $a\leq c \leq d$. If $c \notin B$, then there exists $b \in B$ such that $a\leq b\leq c\leq d$. Hence, $A \ind _B D$. Similarly, we also have that if $A \ind _B C$, $D \ind _{AB} C$, then $D \ind _B C $. Other axioms can also be checked easily.
\end{ex}

Throughout this section, we will assume that $\M$ is a homogeneous structure with a stationary weak independence relation.

\begin{dfn}
We say that $g \in Aut(\M)$ \emph{moves almost $R$-maximally} if for any finite set $X$ and $n$-type $p$ over $X$, there is a realisation $a$ of $p$ such that 
\[ a \ind _{X} g(a).\]

We say $g \in Aut(\M)$ \emph{moves almost $L$-maximally} if for any finite set $X$ and $n$-type $p$ over $X$, there is a realisation $a$ of $p$ such that 
\[g(a)  \ind _{X} a.\]
\end{dfn}

We now prove Theorem \ref{main} using a similar approach to \cite{tentziegler2013isometry}.


The following two lemmas generalise Lemma 3.1 in \cite{tentziegler2013isometry} and the proof follows that of \cite{tentziegler2013isometry}.

\begin{lem}\label{addset}
For finite substructures $A,B,C, D \subseteq \M$, if $A \ind _B C$ and $D$ is arbitrary, then 

I. there is some $D' \subseteq \M$ such that $tp(D' / BC )=tp(D/ BC)$ and $A \ind _B CD'$ 

II. there is some $D''\subseteq \M$ such that $tp(D'' / AB )=tp(D/ AB)$ and $AD'' \ind _B C$.
\end{lem}

\begin{proof}
By Existence, there exists $D'$ realising $tp(D/BC)$ such that $A \ind _{BC} D'$. Then by transitivity on $A \ind _B C$ and $A \ind _{BC} D'$, $A \ind _B CD'$.

Similarly, there exists $D'$ realising $tp(D/AB)$ such that $D' \ind _{AB} C$. Then by transitivity, $AD' \ind _B C$.
\end{proof}

\begin{lem}\label{tz31}
Let $A,B,C$ be finite substructure of $\M$ such that $ A \ind _B C$ and $g_1,...,g_n \in G$. Then 

I. there is $e \in G_{(BC)}$ with $A \ind _B C g_1^{e} (C) \cdots g_n^{e} (C)$.

II. there is $f \in G_{(AB)}$ with $A g_1^{f} (A) \cdots g_n^{f} (A) \ind _B C $.
\end{lem}

\begin{proof}
By the previous lemma, there exist $C_1,...,C_n$ such that 
\[ tp(C_1...C_n/BC)=tp(g_1(C)...g_n(C)/BC) \]
and $A \ind  _B CC_1...C_n$. Then there exists  $e \in G_{(BC)}$ such that $C_i=e( g_i(C) )$ and thus, $g_i ^e(C)=C_i$.

Similarly, there exist $A_1,...,A_n$ such that 
\[ tp(A_1...A_n/AB)=tp(g_1(A)...g_n(A)/AB) \]
and $A A_1...A_n \ind  _B C $. Then there exists  $f \in G_{(AB)}$ such that $A_i=f (g_i(A))$.

\end{proof}

The following proposition uses a similar approach as Proposition 3.2 in \cite{tentziegler2013isometry}. Instead of $Y_2 \ind _{Y_3} Y_4$ as in \cite{tentziegler2013isometry}, we require $Y_4 \ind _{Y_3} Y_2$ here. Therefore, only step 4 and 5 differ from the original proof.

\begin{prop}\label{tz32}
Let $g_1,...,g_4 \in Aut(\M)$ and $X_0,...,X_4 \subset \M$ be such that $g_i (X_{i-1})=X_i$. Then for $i=1,...,4$, there are $a_i \in G_{(X_{i-1}X_i)}$ and extensions $X_i \subset Y_i$ such that $g_i ^{a_i}(Y_{i-1})=Y_i$ and $Y_0 \ind _{Y_1} Y_2$, $Y_4 \ind _{Y_3} Y_2$.
\end{prop}

\begin{proof}

Step 1. 
Choose a finite extension $X'_1$ of $X_1$ such that $X_0 \ind _{X'_1} X_2 X_3 X_4$. Such an extension exists as we can choose $X'_1 =X_0 \cup ... \cup X_4$.

Step 2. 
Applying the previous lemma to $A=X_0, B= X'_1, C=X'_1X_2X_3X_4$ and automorphisms $g_2,g_3g_2,g_4g_3g_2$, we obtain $e\in G_{(X'_1X_2X_3X_4)}$ such that 
\[ X_0 \ind_{X'_1} X'_1X_2X_3X_4 g^e_2(X'_1X_2X_3X_4) (g_3g_2)^e (X'_1X_2X_3X_4) (g_4g_3g_2)^e (X'_1X_2X_3X_4).  \]
By monotonicity, we have 
\[ X_0 \ind_{X'_1}  g^e_2(X'_1) (g_3g_2)^e (X'_1) (g_4g_3g_2)^e (X'_1).  \]
Let $X'_2= g^e_2 (X'_1)$, $X'_3= (g_3 g_2) ^e (X'_1)$, $X'_4= (g_4 g_3 g_2) ^e (X'_1)$. Since $X_1 \subseteq X'_1$ and $e\in G_{(X'_1X_2X_3X_4)}$, we have
\[  X'_2= g^e_2 (X'_1)\supseteq g^e_2(X_1)=eg_2(X_1) =e(X_2)=X_2,
\]
\[   X'_3= (g_3 g_2) ^e (X'_1)=g_3 ^e (X'_2) \supseteq g_3 ^e (X_2)=X_3,
\]
\[   X'_4= (g_4g_3 g_2) ^e (X'_1)=g_4 ^e (X'_4) \supseteq g_4 ^e (X_3)=X_4.
\]
Hence, we have 
\[  X_0 \ind_{X'_1}  X'_2 X'_3 X'_4.
\]

Step 3. Applying the previous lemma to $A=X_0 X'_1, B= X'_1, C= X'_2 X'_3 X'_4$ and automorphism $g^{-1}_1$, we have $f\in G_{(X_0X'_1)}$ such that 
\[ X_0 X'_1 (g^{-1})^f( X_0 X'_1) \ind_{X'_1}  X'_2 X'_3 X'_4 .
 \]
By monotonicity, 
\[(g_1^{-1})^f( X'_1) \ind_{X'_1}  X'_2 X'_3 X'_4 .
\]
Let $X'_0=(g_1^{-1})^f( X'_1)$. Since $X_1 \subseteq X'_1$ and $f\in G_{(X_0X'_1)}$, 
\[ X'_0=(g_1^{-1})^f( X'_1) \supseteq (g^{-1}_1)^f( X_1) fg^{-1}(X_1) =f(X_0)=X_0.
\]
Hence we have $X'_0  \ind_{X'_1}X'_2 X'_3 X'_4$.

Set $h_1=g^f_1, h_2=g^e_2,h_3=g^e_3,h_4=g^e_4$, then $h^{-1}_i(X'_{i})=X'_{i-1}$.

Step 4. Extend $X'_3$ to $Y_3$ such that $X'_2 X'_1 \ind _{Y_3} X'_4$.

Applying the previous lemma to $A=X'_4, B=Y_3, C=Y_3X'_2X'_1$ and automorphisms $h^{-1}_3, h^{-1}_2 h^{-1}_3$, we get $ a \in G_{(Y_3 X'_2 X'_1)}$ such that 
\[ X'_4   \ind _{Y_3} Y_3X'_2 X'_1( h_3^{-1})^a (Y_3X'_2 X'_1) (h_2^{-1} h_3^{-1})^a (Y_3X'_2 X'_1).\]

By monotonicity, \[  X'_4  \ind _{Y_3} ( h_3^{-1})^a (Y_3) (h_2^{-1} h_3^{-1})^a (Y_3) . \]

Let $Y_2 =h_3^{-a} (Y_3) \supseteq X'_2$ and $Y_1 =(h_2 ^{-1} h_3^{-1})^a (Y_3) \supseteq X'_1$. Then \[ X'_4 Y_3   \ind _{Y_3} Y_2 Y_1 . \]

Step 5. Applying the previous lemma to $A= X'_4Y_3, B= Y_3, C=Y_1Y_2$ and $h_4$, we have $b \in G_{(Y_3X'_4)}$ such that \[  Y_3X'_4 h_4^b (Y_3X'_4) \ind _{Y_3} Y_1Y_2. \]

By monotonicity, \[ h_4^b (Y_3) \ind _{Y_3}   Y_1Y_2. \]

Let $Y_4 =h_4^b (Y_3) \supseteq X'_4$, then \[ Y_4 \ind _{Y_3} Y_1 Y_2. \]

Step 6. By Lemma \ref{addset}, since $X'_0 \ind _{X'_1} X'_2 X'_3 X'_4$, we can find $Y'_1Y'_2Y'_3Y'_4$ realising $tp(Y_1 Y_2 Y_3 Y_4/ X'_1 X'_2 X'_3 X'_4)$ such that $X'_0 \ind _{X'_1} Y'_1Y'_2 Y'_3 Y'_4$. Then by invariance, we also have $Y'_4 \ind _{Y'_3} Y'_1 Y'_2$. Therefore, we may assume that \[X'_0 \ind_{X'_1} Y_1Y_2Y_3Y_4.\] 

Then, by Monotonicity,
\[X'_0 \ind_{Y_1} Y_2Y_3Y_4.\]

Applying the previous lemma to $A=X_0 Y_1, B= Y_1, C= Y_2 Y_3 Y_4$ and automorphism $h^{-1}_1$, we find some $c\in
  G_{(X'_0Y_1)}$ such that
\[X_0 Y_1 h_1 ^{-c} (X_0 Y_1) \ind_{Y_1} Y_2Y_3Y_4.\]

Let $Y_0=h_1^{-c}(Y_1)$. By Monotonicity, \[Y_0 \ind_{Y_1} Y_2Y_3Y_4.\]

Altogether, we have $g_1 ^{cf} (Y_0)=Y_1, g_2 ^{ae} (Y_1)=Y_2, g_3 ^{ae}(Y_2)=Y_3, g_4 ^{be}(Y_3)=Y_4$.

\end{proof}

Lemma \ref{tz35} follows a similar approach as Lemma 3.5 in \cite{tentziegler2013isometry}. Lemma \ref{tz36} requires some modifications from Lemma 3.6 in \cite{tentziegler2013isometry}. Combining these two lemmas, we obtain Proposition \ref{tz34} the same way as obtaining Proposition 3.4 from Lemma 3.5 and 3.6 in \cite{tentziegler2013isometry}.

\begin{lem}\label{tz35}
Let $g\in G$ move almost $L$-maximally, let $X,Y,C$ be finite sets such that
  $g(X)=Y$ and $X \ind _Y C$ and let $x$ be a tuple. Then there is some
  $a\in G_{(XY)}$  such that
  \[ g^{a} (x) \ind _Y C .\]
\end{lem}
\begin{proof}
Let $x'$ be a realisation of $tp(x/XY)$ moved almost $L$-maximally by $g$ and
  let $a_1\in  G_{(XY)}$ be such that $a_1(x')=x$. So we have $g (x' )\ind _{XY} x'$, which is the same as $ ga^{-1}_1 (x)\ind _{XY} a^{-1}_1 (x)$. Acting by $a_1$ on it, we have
\[ g^{a_1} (x) \ind _{ XY} x. \]

By existence, there exists $y$ realising $tp(g^{a_1}(x)/XYx)$ with \[ y \ind_{xXY} C.\] 
By invariance, we have $y \ind _{XY} x$.

By transitivity on $y \ind _{XY} x$ and $y \ind _{xXY} C$, we have \[y \ind _{XY} C.\] Together with $X \ind _Y C$, we have by transitivity, \[ y \ind _Y C.\]

Let $a_2\in G_{(xXY)}$ with $a_2g^{a_1}(x)=g^{a_2a_1}(x)=y$.

\end{proof}

\begin{rmk}
We can also prove the following by using a symmetric argument:

Let $g\in G$ move almost $R$-maximally, let $X,Y,C$ be finite sets such that
  $g(X)=Y$ and $C \ind _Y X$ and let $x$ be a tuple. Then there is some
  $a\in G_{(XY)}$  such that
  \[C \ind _Y g^a(x) .\]
\end{rmk}

\begin{lem}\label{tz36}
 Let $g\in G$ moves almost $R$-maximally and let $X,Y$ be finite sets such that
  $g(X)=Y$. Let $x$ and $y$ be tuples satisfying $g(tp(x/X))=tp(y/Y)$ and $x \ind _X yY$, $y \ind _Y X$.  Then
  there is some $a\in G_{(XY)}$ such that
 \[g^a(x)=y.\]
\end{lem}
\begin{proof}
By existence, there exists $y'$ realising $tp(y/XY)$ such that 
\begin{align}\label{361}
y' \ind _{XY} g(Y) .
\end{align}

Since $g$ moves almost $R$-maximally, there exists $y''$ realising $tp(y'/ XYg(Y))$ such that
\begin{align}\label{almostmax}
y'' \ind _{XYg(Y)} g(y'') .
\end{align}

By invariance, we have $y'' \ind _{XY} g(Y) $. By transitivity, we have
\begin{align}\label{363}
y'' \ind _{XY} g(y'') g(Y).
\end{align} 

By invariance, we also have $y'' \ind _Y X$. Then by transitivity on it and (\ref{363}),
\[
y'' \ind _{Y} g(y'') g(Y) .
\]

Then we have $g^{-1}(y'' )\ind _{X} y''Y$.

Since $tp(y''/XY)=tp(y'/XY)=tp(y/XY)$, there exists $b \in G_{(XY)}$ such that $b(y)=y''$, so we have 
\[
g^{-b}(y'' )\ind _{X} yY.
\]

We also have $tp(g^{-b}(y)/X)=g^{-1}tp(y/Y)=tp(x/X)$. By stationarity on $x\ind _{X} yY$ and $g^{-b}(y'' )\ind _{X} yY$, we obtain $tp(g^{-b}(y)/XYy)=tp(x/XYy)$. Therefore, there exists $c \in G_{(XYy)}$ such that $x=c \cdot g^{-b}(y)=g^{-cb} (y)$.
\end{proof}

\begin{prop}\label{tz34}
  Let $g_1,\ldots,g_4\in G$ be such that $g_1,g^{-1}_4$ move almost $L$-maximally and $g^{-1}_2,g_3$ move almost $R$-maximally. Let $Y_0,\ldots,Y_4$ be
  finite sets such that $g_i(Y_{i-1})=Y_i$ for $i=1,\ldots 4$. Assume
  also that $Y_0 \ind_{Y_1}Y_2$ and $Y_4 \ind_{Y_3} Y_2$. Let
  $x_0$ and $x_4$ be two tuples such that $g_4g_3g_2g_1$ maps
  $tp(x_0/Y_0)$ to $tp(x_4/Y_4)$. Then for $i=1,\ldots 4$, there are
  $a_i\in G_{(Y_{i-1}Y_i)}$ such that
  \[g_4^{a_4}\ldots g_1^{a_1}(x_0)=x_4.\]
\end{prop}
\begin{proof}
 Since $g_1$ and $g_4 ^{-1}$ move almost $L$-maximally, by Lemma \ref{tz35}, there exist $a_0\in G_{(Y_0Y_1)}$ and $a_4\in G_{(Y_3Y_4)}$ such that for $x_1=g_1^{a_1}(x_0), x_3=g_4^{-a_4}(x_4)$, we have
 \[x_1 \ind_{Y_1} Y_2, x_3 \ind_{Y_3} Y_2.\] 
By existence, there exists $x_2$ realising the type
  $g_2(tp(x_1/Y_1))=g_3^{-1}(tp(x_3/Y_3))$ such that
  \[x_2 \ind_{Y_2}  x_1Y_1x_3Y_3 .\] 

By Lemma \ref{tz36}, since $g_2^{-1}$ and $g_3$ move almost $R$-maximally, there exist $a_2\in
  G_{(Y_1Y_2)}$ and $a_3\in G_{(Y_2Y_3)}$ such that
  $g_2^{a_2}(x_1)=x_2=(g_3^{-1})^{a_3}(x_3)$.
\end{proof}

\begin{rmk}
It is shown in Lemma 2.8 in \cite{tentziegler2013isometry} that if a countable structure $\M$ has a stationary independence relation, then $Aut(\M)$ has a dense conjugacy class. The proof does not depend on symmetry, so the same statement holds for a countable structure $\M$ with a stationary weak independence relation.
\end{rmk}

\begin{prop}
Suppose $\M$ is a countable structure with a stationary weak independence relation and $g \in Aut (\M)$ be such that $g$ and $g^{-1}$ move both almost $R$-maximally and almost $L$-maximally. 
Define \[\phi: G^4\rightarrow G, \phi:(h_1,\ldots h_4)\mapsto g^{h_1}\ldots g^{h_4}.\] 
Then for any non-empty open set $U\subseteq G^4$, there is some open set $W\subseteq G$ such that $\phi(U)$ is dense in $W$.  
\end{prop}
\begin{proof}
Let $U \subseteq G^4$ be open and non-empty. Assume that $U=U_1\times\ldots \times U_4$, where each $U_i=\mathcal U(u_i)=\{g\in G | u_i \subset g\}$ is the set of automorphisms that extend some finite partial isomorphism $u_i$. Extend each $u_i$ to some $a_i\in G$ such that $g^{a_i}(X_{i-1})=X_i$ for $i=1,\ldots,4$, where $Im(u_i)\subset X_i$. 

By Proposition \ref{tz32}, there exist $b_i\in G_{(X_{i-1}X_i)}$ and extensions $X_i\subset Y_i$ such that $g^{b_ia_i}(Y_{i-1})=Y_i$ and $Y_0 \ind_{Y_1}Y_2$ and $Y_4 \ind_{Y_3} Y_2$.

Let $w$ be the finite isomorphism $g^{b_4a_4}\ldots g^{b_1a_1}$ restricted to $Y_0$ and $W=\mathcal U(w)$. We want to show that $\phi(U)$ is dense in $W$. For an extension $w\subset w'$, let $x$ be an enumeration of $dom(w')\setminus Y_0$ and $y=w'(x)$. By Proposition \ref{tz34}, we obtain $c_i\in G_{(Y_{i-1}Y_i)}$ such that $g^{c_4b_4a_4}\ldots g^{c_1b_1a_1}(x)=y$. Since $b_i$ and $c_i$ both fix $Im(u_i)$ pointwise, we have $c_ib_ia_i\in U_i$. So the $4$-tuple $(c_1b_1a_1,\ldots,c_4b_4a_4)$ belongs to $U$ and is mapped by $\phi$ to $g^{c_4b_4a_4}\ldots g^{c_1b_1a_1}$, which belongs to $W'$.
\end{proof}

The remaining steps in the proof of Theorem \ref{main} now follow exactly as in \cite{tentziegler2013isometry}. We now prove some lemmas which enable us to produce automorphisms that satisfy the hypothese of Theorem \ref{main}.




\begin{lem}
Let $\M$ be a countable structure with a stationary weak independence relation and $g \in Aut (\M)$ move both almost $R$-maximally and almost $L$-maximally. Then there exists $k \in Aut(\M)$ such that $[g,k]$ and its inverse move both almost $R$-maximally and almost $L$-maximally.
\end{lem}
\begin{proof}
We use a back-and-forth construction. Suppose at some stage, we have a partial isomorphism $\tilde{k}:A \rightarrow B$.  Given some type $p$ over $X$, we may assume $p$ is non-algebraic and $X \cup g(X) \subseteq A$, $\tilde{k} X \subseteq g^{-1}B$.

Step I. We want to extend $\tilde{k}$ so that $[\tilde{k},g]$ moves $p$ almost $R$-maximally.

By the existence axiom of $\ind$, there exists a realisation $a'$ of $p$ such that $a' \ind _X A$. Since $g$ moves almost $R$-maximally, there exists $a$ realising $tp(a/A)$ such that $a \ind_A ga$. Since $tp(a/A)=tp(a'/A)$, we also have $a \ind _X A$. By existence, there exists a realisation $b$ of $\tilde{k} \cdot tp(a/A)$ such that $b \ind _B g^{-1}B$. Extend $\tilde{k}$ by sending $a$ to $b$. By acting $\tilde{k}$ on $a \ind _X A$, we get $b \ind _{\tilde{k} X} B$.

Again by the existence axiom, there exists a realisation $c$ of $ \tilde{k}^{-1} \cdot  tp(gb/bB)$ such that $c \ind _{aA} ga$. Since $a \ind_A ga$, by transitivity, we have $c \ind _A ga$. Extend $\tilde{k}$ by sending $c$ to $gb$. 

By transitivity on $b \ind _{\tilde{k} X} B$ and $b \ind _B g^{-1}B$, we have $b \ind _{\tilde{k} X}g^{-1} B$. Acting by $\tilde{k}^{-1}g$ on it, we then get $\tilde{k}^{-1}gb \ind _{\tilde{k}^{-1}g\tilde{k} X} \tilde{k}^{-1}B$, which can be simplified to $c \ind _{\tilde{k}^{-1}g \tilde{k} X} A$.

Since $ \tilde{k}^{-1}g \tilde{k} X \subseteq A $, we can apply transitivity on $c \ind _{\tilde{k}^{-1}g \tilde{k} X} A$ and $c \ind _A ga$ to obtain $ c \ind _{\tilde{k}^{-1}g\tilde{k} X} ga$. Acting by $\tilde{k}^{-1}g^{-1}  \tilde{k}$ on it, we have the required result, i.e.
\begin{align*}
a \ind _X [\tilde{k},g] a.
\end{align*}

Step II. Under the assumption that $g$ moves $L$-maximally and by the exactly same argument, but swapping the sides of $\ind$ as in Step I, we extend $\tilde{k}$ so that $[\tilde{k},g]$ moves $p$ almost $L$-maximally. 

Step III. We want to extend $\tilde{k}$ so that $[g,\tilde{k}]$ moves $p$ almost $R$-maximally. 

By the existence axiom of $\ind$, there exists a realisation $a$ of $p$ such that $a \ind _X g^{-1}(A)$. Since $g$ moves almost $R$-maximally, there exists $b$ realising $\tilde{k}\cdot tp(a/A)$ such that $b \ind _B gb$. Extend $\tilde{k}$ by sending $a$ to $b$. Again by existence, there exists $c$ realising $\tilde{k} \cdot tp(ga/aA)$ such that $c \ind _{bB} gb$. Extend $\tilde{k}$ by sending $ga$ to $c$.

By transitivity on $b \ind _B gb$ and $c \ind _{bB} gb$, we have $c \ind _{B} gb$. By acting $\tilde{k}g$ on $a \ind _X g^{-1}(A)$, we obtain $c \ind _{\tilde{k} g(X)} B$. Again by transitivity and that $\tilde{k}g(X) \subseteq B$, we have $c\ind _{\tilde{k}g(X)} gb$. By acting $g^{-1} \tilde{k}^{-1}$ on it, we obtain $a \ind _X [g,\tilde{k}]a$.

Step IV. Under the assumption that $g$ moves $L$-maximally and by the exactly same argument as in Step III, but swapping the sides of $\ind$, we extend $\tilde{k}$ so that $[g,\tilde{k}]$ moves $p$ almost $L$-maximally. 

Let $k$ be the union of $\tilde{k}$, we have constructed $k \in Aut(\M)$ such that $[\tilde{k},g]$ and $[\tilde{k},g]^{-1}$ move almost $L$-maximally and $R$-maximally.
\end{proof}

Therefore, by Theorem \ref{main}, we have the following corollary:

\begin{coro}\label{corocherlin}
Suppose $\M$ is a countable structure with a stationary weak independence relation and $g \in Aut (\M)$ move both almost $R$-maximally and almost $L$-maximally. Then any element of $G$ is a product of 16 conjugates of $g$ and $g^{-1}$.
\end{coro}

\begin{lem}
Let $\M$ be a countable structure with a stationary weak independence relation and $g \in Aut (\M)$ is such that $g$ moves almost $R$-maximally and $g^{-1}$moves almost $L$-maximally. Then there exists $k \in Aut(\M)$ such that $[g,k]$ and its inverse move both almost $R$-maximally and almost $L$-maximally.
\end{lem}
\begin{proof}
We use a back-and-forth construction. Suppose at some stage, we have a partial isomorphism $\tilde{k}:A \rightarrow B$.  Given some type $p$ over $X$, we may assume $p$ is non-algebraic and $X \cup g(X) \subseteq A$, $\tilde{k} X \subseteq g^{-1}B$. By Step I. and III. in the proof of previous lemma, we know that if $g$ moves almost $R$-maximally, we can extend $\tilde{k}$ so that $[\tilde{k},g]$ and $[\tilde{k},g]^{-1}$ move $p$ almost $R$-maximally. Therefore, we only have to extend $\tilde{k}$ so that $[\tilde{k},g]$ and $[\tilde{k},g]^{-1}$ move $p$ almost $L$-maximally.

Now we extend $\tilde{k}$ so that $[\tilde{k},g]$ moves $p$ almost $L$-maximally.

By existence, there exists $a$ realising $p$ such that $A \ind _X a$ and $b'$ realising $\tilde{k}\cdot tp(a/A)$ such that $g^{-1}(B) \ind_B b'$. Since $g^{-1}$ moves almost $L$-maximally, there exists $b$ realising $tp(b/Bg^{-1}(B))$ such that $g^{-1}(b) \ind _{Bg^{-1}(B)} b$. By invariance and transitivity, we have $g^{-1}(bB) \ind _B b$. Extend $\tilde{k}$ by sending $a$ to $b$. By existence, there exists $c$ realising $\tilde{k}^{-1} \cdot tp(gb/bB)$ such that $ga \ind _{aA} c$. Extend $\tilde{k}$ by sending $c$ to $gb$. 

Then we obtain $B \ind _{\tilde{k}(X)} b$ by acting by $\tilde{k}$ on $A \ind _X a$. By transitivity and $\tilde{k}(X) \subseteq B$, we have $g^{-1}(bB) \ind _{\tilde{k}(X)} b$. Acting by $\tilde{k}^{-1}g$ on it, we obtain $aA \ind_{\tilde{k}^{-1}g\tilde{k}(X)} c$. By transitivity on it and $ga \ind _{aA} c$, we have $ga \ind_{\tilde{k}^{-1}g\tilde{k}(X)} c$. Acting by $\tilde{k}^{-1}g^{-1}\tilde{k}$ on it, we get $[\tilde{k},g] a \ind _X a$.

Now we extend $\tilde{k}$ so that $[g,\tilde{k}]$ moves $p$ almost $L$-maximally.

By existence, there exists $a'$ realising $p$ such that $g^{-1}(A) \ind _X a'$. Since $g^{-1}$ moves $L$-maximally, there exists $a$ realising $tp(a'/g^{-1}(A)X)$ such that $g^{-1}(a) \ind _{g^{-1}(A)X} a$. Then by invariance and transitivity, we have $g^{-1}(a) \ind _X a$ and hence, $a\ind_X ga$. By existence, we can find $b$ realising $\tilde{k} tp(a/A)$ and $c$ realising $\tilde{k} ^{-1}tp(gb/bB)$ such that $c \ind_{aA}ga$. Since $g(X) \subseteq A$, by transitivity, we have $c\ind_{g(X)} ga$. Therefore, acting by $g^{-1}$ on it, we obtain $[g,\tilde{k}] a \ind _X a$.

Let $k$ be the union of $\tilde{k}$, we have constructed $k \in Aut(\M)$ such that $[\tilde{k},g]$ and $[\tilde{k},g]^{-1}$ move almost $L$-maximally and $R$-maximally.
\end{proof}

Therefore, we have the following corollary:

\begin{coro}\label{corodlo}
Suppose $\M$ is a countable structure with a stationary weak independence relation and $g \in Aut (\M)$ is such that $g$ moves almost $R$-maximally and $g^{-1}$moves almost $L$-maximally. Then any element of $G$ is a product of 16 conjugates of $g$ and $g^{-1}$.
\end{coro}

\section{Cherlin's Examples} 

In the appendix of \cite{cherlin1998classification}, Cherlin introduced some special amalgamation classes, called \emph{semi-free amalgamation classes}, defined as follows:

\begin{dfn}\label{semifree}
Given a relational language $\lan$, let $\cl$ be an amalgamation class of finite $\lan$-structures. We say $\mathcal{C}$ is a \emph{semi-free amalgamation class} if there exists $\lan ' \subsetneq \lan$ such that for any finite structures $A,B,C \in \mathcal{C}$ and embeddings $f_1:B \rightarrow A, f_2:B\rightarrow C$, there exist $D\in \mathcal{C}$ and embeddings $g_1:A \rightarrow D, g_2:C\rightarrow D$ such that $g_1f_1(B)=g_2f_2(B)=g_1(A) \cap g_2(C)$ and for any $a \in g_1(A)\setminus g_1 f_1(B), c\in g_2(C)\setminus g_2 f_2(B)$, if $a,c$ are related by some $R \in \lan$, then $R \in \lan'$. We call $\lan'$ the \emph{set of solutions}.

We say a homogeneous structure is \emph{semi-free} if it is the Fra{\"\i}ss{\'e} limit of a semi-free amalgamation class. We can see that a free amalgamation class is a special case of semi-free amalgamation classes where $|\lan'|=1$.
\end{dfn} 


For Cherlin's examples, we fix $\lan$ to be a language consisting of binary and irreflexive relations. We say an $\lan$-structure $A$ is \emph{complete} if every two distinct elements $a,b \in A$ are related by exactly one relation of $\lan$. We denote this relation by $\rel (a,b)$. We will study amalgamation classes of complete $\lan$-structures. Cherlin's examples are specified by sets of forbidden triangles constraints, defined as the following. 

\begin{dfn}\label{forbs}
An $\lan$-structure is a \emph{triangle} if it is a complete structure on three points. Let $S$ be a set of triangles. We define $Forb_c(S)$ to be the set of all complete structures that do not embed any triangle from $S$. We call $S$ \emph{the set of forbidden triangles} of $Forb_c(S)$.
\end{dfn}

We can think of the structures in $Forb_c(S)$ as complete edge-coloured directed graphs that do not embed some coloured triangles by taking the elements of the structures as vertices and the relations as colours. In this section, we let $S$ be a set of forbidden triangles such that the corresponding $Forb_c(S)$ is a semi-free amalgamation class. Then we can take its Fra{\"\i}ss{\'e} limit and denote it by $\M_S$. 

Cherlin provided some examples of $S$ where $\lan$ consists of three or four symmetric relations and where $\lan$ consists of two aysmmetric relations in the appendix of \cite{cherlin1998classification}. We now look at the asymmetric examples of general type, listed below with the original indexing. Note that Cherlin called these structures \emph{2-multi-tournaments} in \cite{cherlin2019book}. Consider the following lists of forbidden triangles where the language $\lan$ consists of two asymmetric relations $R$ and $G$. Note that each pair of vertices has exactly one direction, i.e. for all $a,b \in M_S$, $R^{+}(a,b)$ if and only if $R^{-}(b,a)$. We set $R^{+}(a,b)$ if there is a directed $R$-edge from $a$ to $b$. By a triangle $abc$ with relation $R^+R^+R^+$, we mean $R^+(a,b)R^+(a,c)R^+(b,c)$.

\vspace{0.3cm}

$\lan =\{ R^{\pm},G^{\pm} \}$

\# 8 $G^{+}G^{-}G^{+}$  $R^{+}G^{-}G^{+}$ 

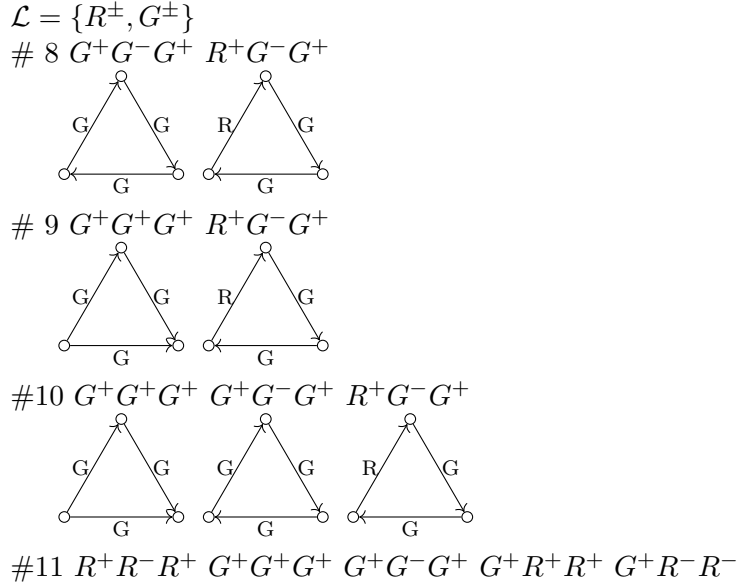
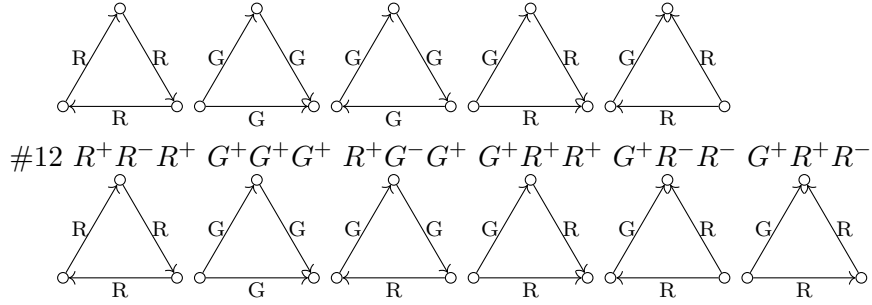
\hspace{0.5cm}
\begin{tikzpicture}
    \tikzstyle{every node}=[draw,circle,fill=white,minimum size=4pt,
                            inner sep=0pt]
  \draw (0,0) node (1) {}
     ++(300:1.5cm) node (2){} 
   ++(180:1.5cm) node (3)  {}
   ++(180:0.4cm) node (4) {}
     ++(120:1.5cm) node (5) {}
   ++(240:1.5cm) node (6)  {};

\draw[->] (1) -- (2) node[draw=none,fill=none,font=\scriptsize,midway,right]{G};
\draw[->]  (2) -- (3)  node[draw=none,fill=none,font=\scriptsize,midway,below] {G};
\draw[->]  (3) -- (1)  node[draw=none,fill=none,font=\scriptsize,midway,left] {R};
\draw[->] (5) -- (4) node[draw=none,fill=none,font=\scriptsize,midway,right]{G};
\draw[->]  (4) -- (6)  node[draw=none,fill=none,font=\scriptsize,midway,below] {G};
\draw[->]  (6) -- (5)  node[draw=none,fill=none,font=\scriptsize,midway,left] {G};

\end{tikzpicture}

\# 9 $G^{+}G^{+}G^{+}$   $R^{+}G^{-}G^{+}$

\hspace{0.5cm}
\begin{tikzpicture}
    \tikzstyle{every node}=[draw,circle,fill=white,minimum size=4pt,
                            inner sep=0pt]
  \draw (0,0) node (1) {}
     ++(300:1.5cm) node (2){} 
   ++(180:1.5cm) node (3)  {}
    ++(180:0.4cm) node (7) {}
     ++(120:1.5cm) node (8) {}
   ++(240:1.5cm) node (9)  {};

\draw[->] (1) -- (2) node[draw=none,fill=none,font=\scriptsize,midway,right]{G};
\draw[->]  (2) -- (3)  node[draw=none,fill=none,font=\scriptsize,midway,below] {G};
\draw[->]  (3) -- (1)  node[draw=none,fill=none,font=\scriptsize,midway,left] {R};
\draw[->] (8) -- (7) node[draw=none,fill=none,font=\scriptsize,midway,right]{G};
\draw[->]  (9) -- (7)  node[draw=none,fill=none,font=\scriptsize,midway,below] {G};
\draw[->]  (9) -- (8)  node[draw=none,fill=none,font=\scriptsize,midway,left] {G};
\end{tikzpicture} 

\#10 $G^{+}G^{+}G^{+}$   $G^{+}G^{-}G^{+}$ $R^{+}G^{-}G^{+}$ 

\hspace{0.5cm}
\begin{tikzpicture}
    \tikzstyle{every node}=[draw,circle,fill=white,minimum size=4pt,
                            inner sep=0pt]
  \draw (0,0) node (1) {}
     ++(300:1.5cm) node (2){} 
   ++(180:1.5cm) node (3)  {}
   ++(180:0.4cm) node (4) {}
     ++(120:1.5cm) node (5) {}
   ++(240:1.5cm) node (6)  {}
    ++(180:0.4cm) node (7) {}
     ++(120:1.5cm) node (8) {}
   ++(240:1.5cm) node (9)  {};

\draw[->] (1) -- (2) node[draw=none,fill=none,font=\scriptsize,midway,right]{G};
\draw[->]  (2) -- (3)  node[draw=none,fill=none,font=\scriptsize,midway,below] {G};
\draw[->]  (3) -- (1)  node[draw=none,fill=none,font=\scriptsize,midway,left] {R};
\draw[->] (5) -- (4) node[draw=none,fill=none,font=\scriptsize,midway,right]{G};
\draw[->]  (4) -- (6)  node[draw=none,fill=none,font=\scriptsize,midway,below] {G};
\draw[->]  (6) -- (5)  node[draw=none,fill=none,font=\scriptsize,midway,left] {G};
\draw[->] (8) -- (7) node[draw=none,fill=none,font=\scriptsize,midway,right]{G};
\draw[->]  (9) -- (7)  node[draw=none,fill=none,font=\scriptsize,midway,below] {G};
\draw[->]  (9) -- (8)  node[draw=none,fill=none,font=\scriptsize,midway,left] {G};
\end{tikzpicture}

\#11 $R^{+}R^{-}R^{+}$ $G^{+}G^{+}G^{+}$  $G^{+}G^{-}G^{+}$  $G^{+}R^{+}R^{+}$ $G^{+}R^{-}R^{-}$ 

\hspace{0.5cm}
\begin{tikzpicture}
    \tikzstyle{every node}=[draw,circle,fill=white,minimum size=4pt,
                            inner sep=0pt]
  \draw (0,0) node (1) {}
     ++(300:1.5cm) node (2){} 
   ++(180:1.5cm) node (3)  {}
   ++(180:0.3cm) node (4) {}
     ++(120:1.5cm) node (5) {}
   ++(240:1.5cm) node (6)  {}
    ++(180:0.3cm) node (7) {}
     ++(120:1.5cm) node (8) {}
   ++(240:1.5cm) node (9)  {}
    ++(180:0.3cm) node (10) {}
     ++(120:1.5cm) node (11) {}
   ++(240:1.5cm) node (12)  {}
   ++(0: 7.2cm) node(13){}
    ++(60:1.5cm) node(14){}
   ++(300:1.5cm) node(15){};

\draw[->] (1) -- (2) node[draw=none,fill=none,font=\scriptsize,midway,right]{R};
\draw[->]  (3) -- (2)  node[draw=none,fill=none,font=\scriptsize,midway,below] {R};
\draw[->]  (3) -- (1)  node[draw=none,fill=none,font=\scriptsize,midway,left] {G};
\draw[->] (5) -- (4) node[draw=none,fill=none,font=\scriptsize,midway,right]{G};
\draw[->]  (4) -- (6)  node[draw=none,fill=none,font=\scriptsize,midway,below] {G};
\draw[->]  (6) -- (5)  node[draw=none,fill=none,font=\scriptsize,midway,left] {G};
\draw[->] (8) -- (7) node[draw=none,fill=none,font=\scriptsize,midway,right]{G};
\draw[->]  (9) -- (7)  node[draw=none,fill=none,font=\scriptsize,midway,below] {G};
\draw[->]  (9) -- (8)  node[draw=none,fill=none,font=\scriptsize,midway,left] {G};
\draw[->] (11) -- (10) node[draw=none,fill=none,font=\scriptsize,midway,right]{R};
\draw[->]  (10) -- (12)  node[draw=none,fill=none,font=\scriptsize,midway,below] {R};
\draw[->]  (12) -- (11)  node[draw=none,fill=none,font=\scriptsize,midway,left] {R};
\draw[->] (15) -- (14) node[draw=none,fill=none,font=\scriptsize,midway,right]{R};
\draw[->]  (15) -- (13)  node[draw=none,fill=none,font=\scriptsize,midway,below] {R};
\draw[->]  (13) -- (14)  node[draw=none,fill=none,font=\scriptsize,midway,left] {G};
\end{tikzpicture}

\#12 $R^{+}R^{-}R^{+}$ $G^{+}G^{+}G^{+}$  $R^{+}G^{-}G^{+}$  $G^{+}R^{+}R^{+}$ $G^{+}R^{-}R^{-}$ $G^{+}R^{+}R^{-}$ 

\hspace{0.5cm}
\begin{tikzpicture}
    \tikzstyle{every node}=[draw,circle,fill=white,minimum size=4pt,
                            inner sep=0pt]
  \draw (0,0) node (1) {}
     ++(300:1.5cm) node (2){} 
   ++(180:1.5cm) node (3)  {}
   ++(180:0.3cm) node (4) {}
     ++(120:1.5cm) node (5) {}
   ++(240:1.5cm) node (6)  {}
    ++(180:0.3cm) node (7) {}
     ++(120:1.5cm) node (8) {}
   ++(240:1.5cm) node (9)  {}
    ++(180:0.3cm) node (10) {}
     ++(120:1.5cm) node (11) {}
   ++(240:1.5cm) node (12)  {}
   ++(0: 7.2cm) node(13){}
    ++(60:1.5cm) node(14){}
   ++(300:1.5cm) node(15){}
  ++(0: 0.3cm) node(16){}
    ++(60:1.5cm) node(17){}
   ++(300:1.5cm) node(18){};

\draw[->] (1) -- (2) node[draw=none,fill=none,font=\scriptsize,midway,right]{R};
\draw[->]  (3) -- (2)  node[draw=none,fill=none,font=\scriptsize,midway,below] {R};
\draw[->]  (3) -- (1)  node[draw=none,fill=none,font=\scriptsize,midway,left] {G};
\draw[->] (5) -- (4) node[draw=none,fill=none,font=\scriptsize,midway,right]{G};
\draw[->]  (4) -- (6)  node[draw=none,fill=none,font=\scriptsize,midway,below] {R};
\draw[->]  (6) -- (5)  node[draw=none,fill=none,font=\scriptsize,midway,left] {G};
\draw[->] (8) -- (7) node[draw=none,fill=none,font=\scriptsize,midway,right]{G};
\draw[->]  (9) -- (7)  node[draw=none,fill=none,font=\scriptsize,midway,below] {G};
\draw[->]  (9) -- (8)  node[draw=none,fill=none,font=\scriptsize,midway,left] {G};
\draw[->] (11) -- (10) node[draw=none,fill=none,font=\scriptsize,midway,right]{R};
\draw[->]  (10) -- (12)  node[draw=none,fill=none,font=\scriptsize,midway,below] {R};
\draw[->]  (12) -- (11)  node[draw=none,fill=none,font=\scriptsize,midway,left] {R};
\draw[->] (15) -- (14) node[draw=none,fill=none,font=\scriptsize,midway,right]{R};
\draw[->]  (15) -- (13)  node[draw=none,fill=none,font=\scriptsize,midway,below] {R};
\draw[->]  (13) -- (14)  node[draw=none,fill=none,font=\scriptsize,midway,left] {G};
\draw[->] (18) -- (17) node[draw=none,fill=none,font=\scriptsize,midway,right]{R};
\draw[->]  (16) -- (18)  node[draw=none,fill=none,font=\scriptsize,midway,below] {R};
\draw[->]  (16) -- (17)  node[draw=none,fill=none,font=\scriptsize,midway,left] {G};
\end{tikzpicture}

\vspace{0.2cm}

Similarly as in \cite{li2018simplicity}, in order to find a stationary weak independence relation and apply Theorem \ref{main} to prove the simplicity of $Aut(\M_S)$, we define the following notion. Note that in order to define the completion process of an amalgamation problem, we say $\lan$ consist of four relations $R^+,R^-,G^+,G^-$ instead of just $R,G$.

\begin{dfn}\label{priority}
Let $\lan$ be a language consisting of $n$ binary and irreflexive relations. Let $ \lan' \subsetneq \{ R^+,R^-|R \in \lan \}$ and suppose $\lan'=\{ R_1,...,R_m\}$. Suppose $\lan'$ is ordered as $R_1 >\cdots >R_m$. For every $A,B,C \in Forb_c(S)$, where $B\subseteq A,C$, define the following way to amalgamate $A$ and $C$ over $B$: for each $a\in A \setminus B, c \in C \setminus B$, first check whether $abc$ form a forbidden triangle for any $b \in B$ if $(a,c) \in R_1$. If $B=\emptyset$ or colouring $(a,c)$ by $R_1$ does not form any forbidden triangle, we let $\rel (a,c) = R_1$. Otherwise, we check the same thing for $(a,c) \in R_2$ and so on. In other words, in the amalgamation, we let $\rel (a,c)= R_i$ where $i$ is the smallest possible integer such that $\rel (a,b) R_i\rel (b,c)  \notin S$ for any $b \in B$.

Denote the resulting amalgamation by $A \otimes _B C$. If $A\otimes _B C$ does not embed any forbidden triangle, i.e. $A\otimes _B C \in Forb_c(S)$, we call it the \emph{prioritised semi-free amalgamation} of $A, C$ over $B$. If for any $A,B,C \in Forb_c(S)$ where $B \subseteq A,C$, $A \otimes _B C \in Forb_c(S)$, then we say $Forb_c(S)$ is a \emph{prioritised semi-free amalgamation class} with respect to the given ordering on $\lan'$.
\end{dfn}

\begin{ex}
For $S$ as in \#8 to \# 10, $Forb_c(S)$ forms a semi-free amalgamation class with $\lan'= \{ R^{+},R^{-} \}$ and $R^+>R^-$. This is because there is no forbidden triangle with two $R$-edges and thus, $AB\otimes _B BC \in Forb_c(S)$ for every $A,B,C \in Forb_c(S)$. 

Taking \# 8 as an example and consider the following amalgamation. Suppose we want to amalgamate $a_1a_2b$ and $bc$ over $b$. For $(a_1,c)$, since completing it with $R^+$ creates a forbidden triangle $R^+G^-G^+$, we complete it with $R^-$. For $(a_2,c)$, completing it with $R^+$ does not create a forbidden triangle, so we let $\rel(a_2,c)=R^+$.

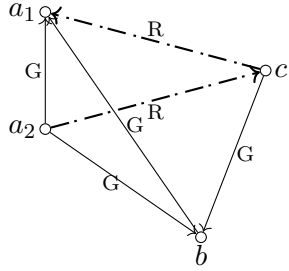
\begin{center}
 \begin{tikzpicture}
    \tikzstyle{every node}=[draw,circle,fill=white,minimum size=4pt,
                            inner sep=0pt]
  \draw (0,0) node (1) [label=left:$a_1$] {}
     ++(345:3cm) node (2) [label=right:$c$]{} 
   ++(195:3cm) node (4) [label=left:$a_2$] {}
   ++(325:2.5cm) node (3)  [label=below:$b$] {};

\draw [->,thick,dash pattern={on 7pt off 2pt on 1pt off 3pt}] (4) -- (2) node[draw=none,fill=none,font=\scriptsize,midway,below] {R};
\draw [->,thick,dash pattern={on 7pt off 2pt on 1pt off 3pt}] (2) -- (1) node[draw=none,fill=none,font=\scriptsize,midway,above] {R};
\draw[->]  (3) -- (1)  node[draw=none,fill=none,font=\scriptsize,midway,right] {G};
\draw[->]  (2) -- (3)  node[draw=none,fill=none,font=\scriptsize,midway,right] {G};
\draw[->]  (4) -- (3)  node[draw=none,fill=none,font=\scriptsize,midway,left] {G};
\draw[->]  (4) -- (1)  node[draw=none,fill=none,font=\scriptsize,midway,left] {G};
\end{tikzpicture}
\end{center}

However, for $S$ as \# 11 and \# 12, $Forb_c(S)$ does not form a prioritised semi-free amalgamation class. As shown in the following examples, we cannot prioritise $R$ or $G$, as otherwise we would have a forbidden triangle while amalgamating $a_1a_2b$ and $bc$ over $b$.

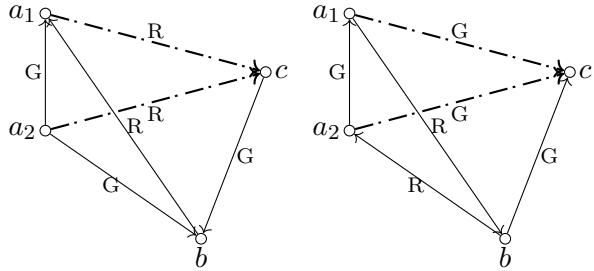
\begin{center}
 \begin{tikzpicture}
    \tikzstyle{every node}=[draw,circle,fill=white,minimum size=4pt,
                            inner sep=0pt]
  \draw (0,0) node (1) [label=left:$a_1$] {}
     ++(345:3cm) node (2) [label=right:$c$]{} 
   ++(195:3cm) node (4) [label=left:$a_2$] {}
   ++(325:2.5cm) node (3)  [label=below:$b$] {}
    ++(0:4cm) node (7) [label=below:$b$] {}
     ++(145:2.5cm) node (8) [label=left:$a_2$]{} 
   ++(15:3cm) node (6) [label=right:$c$] {}
   ++(165:3cm) node (5)  [label=left:$a_1$] {};

\draw [->,thick,dash pattern={on 7pt off 2pt on 1pt off 3pt}] (4) -- (2) node[draw=none,fill=none,font=\scriptsize,midway,below] {R};
\draw [->,thick,dash pattern={on 7pt off 2pt on 1pt off 3pt}] (1) -- (2) node[draw=none,fill=none,font=\scriptsize,midway,above] {R};
\draw[->]  (3) -- (1)  node[draw=none,fill=none,font=\scriptsize,midway,right] {R};
\draw[->]  (2) -- (3)  node[draw=none,fill=none,font=\scriptsize,midway,right] {G};
\draw[->]  (4) -- (3)  node[draw=none,fill=none,font=\scriptsize,midway,left] {G};
\draw[->]  (4) -- (1)  node[draw=none,fill=none,font=\scriptsize,midway,left] {G};

\draw [->,thick,dash pattern={on 7pt off 2pt on 1pt off 3pt}] (8) -- (6) node[draw=none,fill=none,font=\scriptsize,midway,below] {G};
\draw [->,thick,dash pattern={on 7pt off 2pt on 1pt off 3pt}] (5) -- (6) node[draw=none,fill=none,font=\scriptsize,midway,above] {G};
\draw[->]  (5) -- (7)  node[draw=none,fill=none,font=\scriptsize,midway,right] {R};
\draw[->]  (7) -- (6)  node[draw=none,fill=none,font=\scriptsize,midway,right] {G};
\draw[->]  (7) -- (8)  node[draw=none,fill=none,font=\scriptsize,midway,left] {R};
\draw[->]  (8) -- (5)  node[draw=none,fill=none,font=\scriptsize,midway,left] {G};
\end{tikzpicture}
\end{center}

\end{ex}

It was shown in \cite{li2018simplicity} that in the case of symmetric relations, if $Forb_c(S)$ forms a prioritised semi-free amalgamation class, then $\M_S$ has a stationary independence relation. The proof does not depend on symmetry except the proof of the symmetry axiom. Hence, we can generalise the result in \cite{li2018simplicity} to the following theorem and can then apply the results in Sections 2 to the asymmetric cases in \cite{cherlin1998classification}.

\begin{thm}\label{swirthm}
Let $S$ be a set of forbidden triangles such that $Forb_c(S)$ is a prioritised semi-free amalgamation class. Let $\M_S$ be the Fra{\"\i}ss{\'e} limit of $Forb_c(S)$. For any finite substructure $ A,B, C$ of $\M_S$, let $A\ind _B C$ if $ABC=AB \otimes _B BC$ where $AB \otimes _B BC$ is the prioritised semi-free amalgamation defined in Definition \ref{priority}. Then $\ind$ is a stationary weak independence relation on $\M_S$.
\end{thm}

In section 4 of \cite{li2018simplicity}, it was shown that if $S$ satisfies the following condition, then for some $g \in Aut(\M_S)$ and some non-algebraic type $p$, if $g$ fixes its set of realisations pointwise, then $g=1$. This is also true for the asymmetric structures as the arguments do not require symmetry. Therefore, any non-trivial $g \in Aut(\M_S)$ moves infinitely many realisations of $p$.

\begin{con}\label{maincond}
Under the notation and assumptions of Theorem \ref{swirthm}, assume that
\begin{enumerate}[(i)]
\item $S$ does not contain any triangle involving $R_1R_1$ or $R_1 R_2$.
\item Let $a,b, c \in \M_S$ and $B \subseteq \M_S$ be such that $a \ind_{bB} c$. If $\rel (a,b) \in \lan'$, then we have $a \ind _B c$.
\end{enumerate}
\end{con}

We generalise the argument in section 5 of \cite{li2018simplicity} so that for any non-trivial $g \in Aut(\M_S)$, we can find a product of conjugates of $g$ and $g^{-1}$ that moves all types almost $R$-maximally and almost $L$-maximally. Note that we have to work with all $n$-types for all $n$ here rather than only 1-types as in \cite{li2018simplicity}. This is because when applying Tent and Ziegler's original method on structures with a stationary independence relation, we build a \emph{moving maximally} automorphism from a moving almost maximally automorphism. And it is shown in \cite{tentziegler2013isometry} that if an automorphism moves all 1-type maximally, then it moves all $n$-type maximally for any $n$. However, we do not have this property for an automorphism that moves all 1-type almost maximally. The interested reader can check Lemma 2.4 and Lemma 2.6 in \cite{tentziegler2013isometry} for more details. Note that for any non-algebraic type over some $X$, we may assume its realisation does not intersect with $X$.

\begin{lem}\label{colourrange}
Let $S$ be a set of forbidden triangle such that $Forb_c(S)$ is a semi-free amalgamation class. Let $\M_S$ be its Fra{\"\i}ss{\'e} limit. Then given any non-trivial $g \in Aut(\M_S)$, we can construct $h\in Aut(\M_S)$ such that for any non-algebraic type $p$ over some finite set $X$, there exist infinitely many realisations $a$ of $p$ such that $\rel (a_i, [h,g]a_j) \in \lan '$ where $a_i,a_j$ are elements of $\M_S$ appearing as coordinates in the tuple $a$.
\end{lem}
\begin{proof}
List all non-algebraic types over a finite set as $p_1,p_2,...$. We start with the empty map and use a back-and-forth construction to build $h$. Suppose at the some stage, we have a partial isomorphism $\tilde{h}:A \rightarrow B$ such that for any $p_k \in \{p_1,...,p_{n-1}\}$, there exists a realisation $a$ of $p_k$ such that $\rel (a_i, [\tilde{h},g]a_j) \in \lan '$ for all $a_i,a_j \in a $.

Let $p:=p_n$ be a $n$-type over $X$. We want to extend $\tilde{h}$ such that $p$ has a realisation $a$ such that $\rel (a, [\tilde{h},g]a)  \in \lan '$. We may assume $X \subseteq A$ by extending $\tilde{h}$.

Since $p$ is non-algebraic, $p$ has a realisation $a$ that does not intersect with $A \cup g^{-1} (aA) $ and $\tilde{h} \cdot tp(a/A)$ has a realisation $b$ that does not intersect with $B \cup g^{-1} (bB) $. Extend $\tilde{h}$ by sending $a$ to $b$. It is well-defined since $\tilde{h}\cdot tp(a/A)=tp(b/B)$.

By the extension property, there exists a realisation $c$ of $\tilde{h}^{-1} \cdot tp(gb/bB)$ such that $c$ and $ga$ are semi-freely amalgamated over $aA$. Then $(c_i,ga_j)$ is coloured using relations from $\lan'$, i.e. $\rel (c_i,ga_j) \in \lan'$, for all $a_i,a_j \in a$. 

Extend $\tilde{h}$ by sending $c$ to $gb$. Since $\tilde{h} \cdot tp(c/aA)=tp(gb/bB)$, $\tilde{h}$ is a well-defined partal isomorphism. Then $c=\tilde{h}^{-1} g \tilde{h} a$ and we have, for any $a_i,a_j \in a$,
\[ \rel (a_i,[\tilde{h},g]a_j) = \rel ( a_i, \tilde{h}^{-1} g ^{-1} \tilde{h}g a_j )= \rel (\tilde{h}^{-1} g \tilde{h} a_i, g a_j )= \rel (c_i,ga_j) \in \lan'. \] 

At every alternative step, we can make sure $X \subset B$ by extending $\tilde{h}$. Let $h$ be the union of all $\tilde{h}$ over each step, it is an automorphism since it is well defined and bijective as we made sure every finite subset of $\M$ is contained in both domain and image. 
\end{proof}

\begin{thm}
Let $S$ be a set of forbidden triangles satsifying Condition \ref{maincond} and $Forb_c(S)$ the set of all finite complete $\lan$-structures that does not embed any triangle from $S$. Let $\M_S$ be its Fra{\"\i}ss{\'e} limit. Let $g \in Aut(\M_S)$ be an automorphism of $\M_S$ such that for any non-algebraic type $p$ over some finite set $X$, there exist infinitely many realisations $a$ of $p$ with $\rel (a_i, ga_j) \in \lan '$ for all $a_i,a_j \in a$. Then there exists $k \in Aut(\M_S)$ such that $[g,k]$ moves both almost $R$-maximally and almost $L$-maximally.
\end{thm}

\begin{proof}
We use a back-and-forth construction. Suppose at some stage, we have a partial isomorphism $\tilde{k}:A \rightarrow B$. We want to extend $\tilde{k}$ so that $\tilde{k}$ moves some given type $p$ over some finite set $X$ almost $R$-maximally. We may assume $p$ is non-algebraic and realisation of $p$ does not intersect with $X$. We may also assume $X \cup gX \subseteq A$ and $\tilde{k} X \subseteq g^{-1}B$ by extending $\tilde{k}$, then we have $\tilde{k}g X \subseteq B$.

By existence, there exists a realisation $a$ of $p$ such that $a \ind_X g^{-1}A$.

By the assumption on $g$, there exists a realisation $b$ of $\tilde{k}tp(a/A)$ such that $\rel (b_i, gb_j) \in \lan '$ for all $b_i,b_j \in b$. 

There also exists a realisation $c$ of $ \tilde{k} \cdot  tp(ga/aA)$ such that $c \ind_{bB} gb$. Since $\rel (b_i, gb_j) \in \lan '$ for all $b_i,b_j \in b$, by part (ii) of Condition \ref{maincond}, we obtain $c \ind_{B} gb$.

Extend $ \tilde{k}$ by sending $a$ to $b$ and $ga$ to $c$. Then we have $c \ind_{\tilde{k}g(X)} B$ from $a \ind_X g^{-1}A$ by invariance. And by transitivity on $c \ind_{B} gb$ and $c \ind_{\tilde{k}g(X)} B$, we have $c \ind_{\tilde{k}g(X)} gb$. By invariance, we have 
\[ a \ind _X [g,\tilde{k}]a\]

By the same construction, but swapping sides of $\ind$, we can also extend $ \tilde{k} $ such that there is a realisation $a'$ of $p$ such that $[g,\tilde{k}] a' \ind _X a'$.
\end{proof}

In Section 6.1 of \cite{li2018simplicity}, it was shown that if $S$ satisfies the following condition, then $S$ also satisfies Condition \ref{maincond} and $Forb_c(S)$ forms a prioritised semi-free amalgamation class. The argument does not depend on symmetry. Therefore, we have the following theorem by the arguments in this section. Part (iii) follows from (ii) by Corollary \ref{corocherlin}.

\begin{con}\label{condition1}
Let $S$ be a set of forbidden triangles. Assume that $S$  does not contain any triangle of the form $R_i R_j R'$ where $R_i, R_j \in \lan '$ and $R' \in \lan$.

\end{con}

\begin{thm}\label{generalcase}
Let $S$ be a set of forbidden triangles satsifying Condition \ref{condition1} and $Forb_c(S)$ the set of all finite complete $\lan$-structures that does not embed any triangle from $S$. Then 
\begin{enumerate}[(i)]
\item $Forb_c(S)$ forms a prioritised semi-free amalgamation class. Let $\M_S$ be its Fra{\"\i}ss{\'e} limit.
\item There is a stationary weak independence relation on $\M_S$ and for any non-trivial $g \in Aut(\M_S)$, there exists $h,k \in Aut(\M_S)$ such that $[[h,g],k]$ moves almost $R$-maximally and $L$-maximally.
\item $Aut(\M_S)$ is simple.
\end{enumerate}
\end{thm}

Since for $S$ as in \#8 to \# 10, it satisfies Condition \ref{condition1}, we have the following corollary.

\begin{coro}
Let $S$ be as in \#8 to \# 10 and $Forb_c(S)$ the set of all finite complete $\lan$-structures that does not embed any triangle from $S$. Then $Forb_c(S)$ forms a prioritised semi-free amalgamation class with $\lan'= \{ R^{+},R^{-} \}$ and $ R^{+} > R^{-}$. Let $\M_S$ be its Fra{\"\i}ss{\'e} limit. Then $Aut(\M_S)$ is simple.
\end{coro}


\section{Further Work}

The author of this paper is currently working on applying Theorem \ref{main} to linearly ordered structures. The author has shown that the theorem can be applied to the dense linear order and linearly ordered free homogeneous structures. The result on the dense linear order was first proved in \cite{holland1963lattice} and \cite{lloyd1964lattice}. The theorem provides an alternative way of proving the same result. The author has also applied the main theorem to show that the automorphism groups of linearly ordered free homogeneous structures is simple. The same result was proved recently in \cite{ckt2019simplicity} with a somewhat different approach. The author is currently working on generalising the methods so that it can be applied to more general linear ordered structures, for example, the random permutation.

\printbibliography

\Addresses

\end{document}